\documentclass{commat}

\usepackage{mathrsfs}
\newtheorem{fact}{Fact}

\title{%
    On commutativity of prime rings with skew derivations
    }

\author{%
   Nadeem ur Rehman and Shuliang Huang
    }

\affiliation{
    \address{Nadeem ur Rehman --
    Department of Mathematics,  Aligarh Muslim University, Aligarh-202002 India.
        }
    \email{%
    nu.rehman.mm@amu.ac.in
    }
    \address{Shuliang Huang --
    School of Mathematics and Finance, Chuzhou University,  Chuzhou-239000 China.
        }
    \email{%
    shulianghuang@163.com
    }
    }

\abstract{%
   Let $\mathscr{R}$ be a prime ring of Char$(\mathscr{R}) \neq 2$ and $m\neq 1$ be a positive integer. If $S$ is a nonzero skew derivation with an associated automorphism $\mathscr{T}$ of $\mathscr{R}$ such that $([S([a, b]), [a, b]])^{m} = [S([a, b]), [a, b]]$ for all $a, b \in \mathscr{R}$, then $\mathscr{R}$ is commutative.
    }



\VOLUME{31}
\YEAR{2023}
\NUMBER{1}
\firstpage{169}
\DOI{https://doi.org/10.46298/cm.10319}

\begin{paper} 

\section{Introduction}
In all that follows, unless specifically stated otherwise, $\mathscr{R}$ will be an associative ring, $Z(\mathscr{R})$ the center of $\mathscr{R}$, $\mathscr{Q}$ its Martindale quotient ring and $U$ its Utumi  quotient ring. The center $\mathscr{C}$ of $\mathscr{Q}$ or $U$, called the extended centroid of $\mathscr{R}$, is a field (see \cite{3} for further details). For any $a, b
\in \mathscr{R}$, the symbol $[a,b]$  denotes the Lie product $ab-ba$. Recall that a ring $\mathscr{R}$ is prime if  for any $a, b\in \mathscr{R}$,
$a\mathscr{R}b=(0)$ implies $a=0$ or $b=0$, and is semiprime if  for any $a\in
\mathscr{R}$, $a\mathscr{R}a=(0)$ implies $a=0$.  An additive subgroup $\mathscr{L}$
of $\mathscr{R}$ is said to be a Lie ideal of $\mathscr{R}$ if $[l,r]\in \mathscr{L}$ for all
$l\in \mathscr{L}$ and $r\in \mathscr{R}$. By a derivation  of $\mathscr{R}$, we mean an additive map
$d:\mathscr{R}\longrightarrow \mathscr{R}$ such that $d(ab)=d(a)b+ad(b)$
holds for all $a,b\in \mathscr{R}$.  An additive map $F:\mathscr{R} \longrightarrow \mathscr{R}$
is called a generalized derivation if there exists a derivation $d:\mathscr{R}
\longrightarrow \mathscr{R}$ such that $F(ab)=F(a)b+ad(b)$ holds for all
$a,b\in \mathscr{R}$, and $d$ is called the associated derivation of $F$.  The standard identity $s_{4}$ in  four variables is defined as
follows:$$s_{4}=\sum(-1)^{\tau}X_{\tau(1)}X_{\tau(2)}X_{\tau(3)}X_{\tau(4)}$$
where $(-1)^{\tau}$ is the sign of a permutation $\tau$ of the
symmetric group of degree $4$.

 It is well known that any automorphism of $\mathscr{R}$ can be uniquely extended to an automorphism of $\mathscr{Q}$. An automorphism $\mathscr{T}$ of $\mathscr{R}$ is called $\mathscr{Q}$-inner if there exists an invertible element $\alpha \in \mathscr{Q} $ such that $ \mathscr{T}(a) = \alpha a\alpha^{-1}$  for every $a \in \mathscr{R}$. Otherwise, $\mathscr{T}$  is called $\mathscr{Q}$-outer. Following \cite{7},  an additive map $S : \mathscr{R} \to \mathscr{R}$ is
said to be a skew derivation if there exists an automorphism $\mathscr{T}$ of $\mathscr{R}$  such that $S(ab) =S(a)b+\mathscr{T}(a)S(b)$  holds for every $a,b\in \mathscr{R}$. It is easy to see that if $\mathscr{T} =1_{\mathscr{R}}$, where $1_{\mathscr{R}}$ the identity map on $\mathscr{R}$, then a skew derivation is just a usual derivation.  If  $\mathscr{T} \neq 1_{\mathscr{R}}$, then $\mathscr{T} - 1_{\mathscr{R}}$ is a skew derivation. Given any $b \in \mathscr{Q}$, obviously the map $S : a \in  \mathscr{R}  \to  ba -  \mathscr{T}(a)b$ defines a skew derivation of $\mathscr{R}$, called $\mathscr{Q}$-inner skew derivation. If a skew derivation $S$  is not $\mathscr{Q}$-inner, then it is called $\mathscr{Q}$-outer. Hence the concept of skew derivations unites the notions of derivations and automorphisms, which have been examined many algebraists from diverse points of view (see \cite{CL},  \cite{Lanski} and  \cite{Liu}).

 A classical result of Divinsky \cite{11} states that if $\mathscr{R}$ is a simple Artinian ring,  $\sigma$ a non-identity automorphism such that  $[\sigma(a),a]=0$  for all $a\in \mathscr{R}$, then $\mathscr{R}$ must be commutative. Many authors have recently investigated and demonstrated commutativity of prime and semiprime rings using derivations, automorphisms, skew derivations, and other techniques that satisfy specific polynomial criteria (see  \cite{0}, \cite{DR},  \cite{27},  \cite{28},  \cite{31} and references therein). Carini and De Filippis \cite{4},   showed  if   a 2-torsion free semiprime ring $\mathscr{R}$ admits a nonzero derivation $d$ such that  $[d([a,b]),[a,b]]^{n}=0$ for  all $a,b\in \mathscr{R}$, then there exists a central idempotent element $e\subseteq U$ such that on the direct sum decomposition $eU\bigoplus(1-e)U$, $d$ vanishes identically on $eU$ and the ring $(1-e)U$ is commutative.
In  \cite{13},  Scudo and Ansari studied the identity  $[G(u),u]^{n}=[G(u),u]$ involving a nonzero generalized derivation $G$ on a noncentral Lie ideal of a prime ring $\mathscr{R}$ and they described the structure of $\mathscr{R}$.  Wang \cite{32} obtained that if  $\mathscr{R}$ is a prime ring, $\mathscr{L}$ a non-central Lie ideal of $\mathscr{R}$ such $[\sigma(a),a]^{n}=0$  for all $a\in \mathscr{L}$, and if either $char(\mathscr{R}) > n$ or  $char(\mathscr{R})=0$, then $\mathscr{R}$ satisfies $s_{4}$. Replaced the automorphism $\sigma$ by a skew derivation $d$, it is  proved in \cite{9} the following result:  Let $\mathscr{R}$ be a prime ring of characteristic different from $2$ and $3$, $\mathscr{L}$  a non-central Lie ideal of $\mathscr{R}$, $d$ a nonzero skew derivation of $\mathscr{R}$, $n$ is a fixed positive integer. If $[d(a),a]^{n}=0$ for all $a\in \mathscr{L}$, then $\mathscr{R}$ satisfies $s_{4}$.

 Motivated by  the previous cited results, our aim here is to examine what happens if a prime ring $\mathscr{R}$ admits a nonzero skew derivation $S$ such that
 \[
 ([S([a, b]), [a, b]])^{m} = [S([a, b]), [a, b]] 
 \quad \textup{ for all } a, b\in \mathscr{R}.
 \]

\section{Notation and Preliminaries}

First, we mention  some important well-known facts which are  needed in the proof of our results.

\begin{fact}[{\cite[Lemma 7.1]{2}}] \label{f1}
	Let $V_{D}$ be a vector space over a division ring $D$ with $dim V_{D} \geq 2$ and $\phi \in End(V)$. If $r$ and $\phi r$ are $D$-dependent for every $r \in V$, then there exists $\lambda \in D$  such that $\phi r = \lambda r$ for every $r \in V$.
\end{fact}

\begin{fact}[{\cite[Theorem 1]{6}}] \label{f2}
	Let $\mathscr{R}$ be a prime ring and $I$ be a two-sided ideal of $\mathscr{R}$. Then $I$, $\mathscr{R}$ and $\mathscr{Q}$ satisfy the same generalized polynomial identities (GPIs) with automorphisms.
\end{fact}

\begin{fact}[{\cite[Fact 4]{8}}] \label{f3}
	Let $\mathscr{R}$ be a domain and $\mathscr{T}$ be an automorphism of $\mathscr{R}$ which is outer. If $\mathscr{R}$ satisfies a GPI $\Xi(r_{i}, \mathscr{T}(r_{i}))$, then $\mathscr{R}$ also satisfies the nontrivial GPI $\Xi(r_{i}, s_{i})$, where $r_{i}$ and $s_{i}$ are distinct indeterminates.
\end{fact}

  \begin{lemma}\label{l1}
  	Let $\mathscr{R}$ be a dense subring of the ring of linear transformations of a vector space $V$ over a division ring $D$ and  $m\neq 1$ a positive integer. If $\mathscr{T}: \mathscr{R} \to \mathscr{R}$ is an automorphism of $\mathscr{R}$ and $\vartheta \in \mathscr{R}$ such that
  	$$([\vartheta[a, b] - \mathscr{T}([a, b])\vartheta, [a, b]])^{m}  = [\vartheta[a, b] - \mathscr{T}([a, b])\vartheta, [a, b]],$$
  	for every $a, b \in \mathscr{R}$, then $dim_{D}V = 1$.
\end{lemma}
\begin{proof}
	We have
	$$([\vartheta[a, b] - \mathscr{T}([a, b])\vartheta, [a, b]])^{m}  = [\vartheta[a, b] - \mathscr{T}([a, b])\vartheta, [a, b]],$$
	for every $a, b \in \mathscr{R}$. As $\mathscr{R}$ and $\mathscr{Q}$ satisfy the same GPIs  with automorphisms by  Fact \ref{f2}, and hence it is a GPI for $\mathscr{Q}$. We prove it by contradiction. We assume that $dim_{D}V \geq 2$. There exists a semi-linear automorphism $\Phi\in End(V)$, by \cite[ p.79]{20},  such that $\mathscr{T}(a) = \Phi a\Phi^{-1}$ $\forall a \in \mathscr{Q}$. Hence, $\mathscr{Q}$ satisfies
	$$([\vartheta[a, b]-\Phi[a, b]\Phi^{-1}\vartheta, [a, b]])^{m} = [\vartheta[a, b] - \Phi[a, b]\Phi^{-1}\vartheta, [a, b]].$$
	Suppose that $\Phi u \not \in span_{D}\{u, \Phi^{-1}\vartheta u\}$, then $\{u, \Phi u, \Phi^{-1}\vartheta u\}$ is linearly $D$-independent. By density theorem for $\mathscr{R}$, there exists $a, b \in \mathscr{R}$ such that
	$$\begin{array}{llllll}
	au=0& a\Phi^{-1}\vartheta u = 2u & a\Phi u = u\\
	bu = -u & b\Phi^{-1} \vartheta u=0 & b\Phi u = 0.
	\end{array}$$
	The above relation gives $[a, b]u =0$, $[a, b]\Phi^{-1} \vartheta u = 2u$ and $[a, b]\Phi u = u$. This implies that
	$$
 (2^{m} - 2)u = \left(([\vartheta[a, b] - \Phi[a, b]\Phi^{-1}\vartheta, [a, b]])^{m} - [\vartheta[a, b] - \Phi[a, b]\Phi^{-1}\vartheta, [a, b]]\right) u = 0,$$
 a contradiction.

 Now, we assume that $\Phi u\in Span_{D}\{u, \Phi^{-1}\vartheta u\}$, then $\Phi u = u\zeta +\Phi^{-1}\vartheta u\theta$ for some $\zeta, \theta \in D$. We see that $\theta \neq 0$ otherwise if $\theta = 0$, then we get $\Phi u=u\zeta$ and hence this gives that $u = \Phi^{-1}u\zeta$. Again by density theorem for $\mathscr{R}$, $\exists a, b \in \mathscr{R}$, we have
 $$\begin{array}{lll}
 au = 0& a\Phi^{-1}u = 2u\\
 bu=-u & b\Phi^{-1}u = 0.
 \end{array}$$
 The above expression again gives that a contradiction
 $$(2^{m} \theta^{m} - 2 \theta)u = \left(([\vartheta[a, b] - \Phi[a, b]\Phi^{-1}\vartheta, [a, b]])^{m} - [\vartheta[a, b] - \Phi[a, b]\Phi^{-1}\vartheta, [a, b]]\right) u = 0.$$
 For $u\in V$, the set $\{u, \Phi^{-1}\vartheta u \}$ is $D$-dependent. By Fact \ref{f1},  $\exists \Delta \in D$ such that $\Phi^{-1} \vartheta u=u\Delta$, $\forall u \in V$ and hence we have
 $$\mathscr{T}(a) \vartheta u =(\Phi a\Phi^{-1}) \vartheta u = \Phi au \Delta$$
 and $$(\mathscr{T}(a)\vartheta - \vartheta a)u = \Phi(au\Delta)-\vartheta au = \Phi(\Phi^{-1} \vartheta au) -  \vartheta au = 0.$$
 The last expression forces that $(\mathscr{T}(a)\vartheta -  \vartheta a)V = (0)$  $\forall a \in \mathscr{R}$, and  hence $\mathscr{T}(a)V = (0)$  $\forall a \in \mathscr{R}$ and as $V$ is faithful, it yields that $\mathscr{T}(a) = 0$ $\forall a \in \mathscr{R}$. This is a contradiction.
\end{proof}

\section{Main Results}

\begin{proposition}\label{p1}
	Let $m\neq 1$ be  a positive integer,  $\mathscr{R}$ be a prime ring of char$(\mathscr{R}) \neq 2$ and $\vartheta \in \mathscr{Q}$ such that
	$$([\mathscr{T}([a, b])\vartheta, [a, b]])^{m} = [\mathscr{T}([a, b])\vartheta, [a, b]].$$
	Then $\vartheta \in \mathscr{C}$.
\end{proposition}

\begin{proof}
First we assume that $\mathscr{T}$ is an identity automorphism of $\mathscr{R}$. Then we have that $([[a, b]\vartheta, [a, b]])^{m} = [[a, b]\vartheta, [a, b]]$
is a GPI of $\mathscr{R}$.
On contrary we assume that $\vartheta \not\in \mathscr{C}$. Since the identity $([[a, b]\vartheta, [a, b]])^{m} = [[a, b]\vartheta, [a, b]]$ is satisfied by $\mathscr{Q}$ (Fact  \ref{f2}). As $\vartheta \not \in \mathscr{C}$, then the above identity is an non-trivial GPI for $\mathscr{Q}$. By Martindale's theorem in \cite{24}, $\mathscr{Q}$ is primitive ring which is isomorphic to a dense ring of linear transformations of a vector space $V$ over~$\mathscr{C}$.

Assume that dim$\mathscr{C}(V)=l$, where $1<l\in  \mathbb{Z}^+$. For this situation, we take $\mathscr{Q}=M_{l}(\mathscr{C})$
as a ring of $l \times l$ matrices over the field $\mathscr{C}$ such that $([[a, b]\vartheta, [a, b]])^m = [[a, b]\vartheta, [a, b]]$  for all $a, b\in M_{l}(\mathscr{C})$.

Let $e_{ij}$ be the usual unit matrix with $1$ in $(i, j)$-entry and zero elsewhere. First, we claim  that $\vartheta$ is a diagonal matrix. Say $\vartheta = \sum_{ij} e_{ij}\vartheta_{ij}$, where $\vartheta_{ij} \in  \mathscr{C}$. Choose $a = e_{ij}, b = e_{jj}$. Then by the hypothesis, we have $([e_{ij}\vartheta, e_{ij}])^m = [e_{ij}\vartheta, e_{ij}]$, i.e, $e_{ij}\vartheta_{ij} = 0$ and so $\vartheta_{ji} = 0$, for any $i\neq  j$ and hence $\vartheta$ is a diagonal matrix.

Since $\xi \in Aut_{\mathscr{C}}(\mathscr{Q})$, the expression
$$([[a, b]\xi(\vartheta), [a, b]])^{m} = [[a, b]\xi(\vartheta), [a, b]]$$
is also a GPI for $\mathscr{Q}$, therefore $\xi(\vartheta)$ is also diagonal.  The automorphism, in particular $\xi(\vartheta) = (1 +e_{ij})\vartheta(1 -e_{ij})$, for any $i \neq j$ and say
$\vartheta^{\xi} = \sum_{ij} e_{ij}\vartheta'_{ij}$, where $\vartheta'_{ij} \in \mathscr{C}$. Since $\vartheta'_{ij} = 0$, then  we get $0 =\vartheta'_{ij} = \vartheta_{jj} - \vartheta_{ii}$, by easy computation. So that $\vartheta_{jj} = \vartheta_{ii}$ hold for any $i \neq  j$, and we get a contradiction that $\vartheta \in \mathscr{C}$.

Assume that $dim_{\mathscr{C}}V = \infty$.
\begin{equation}\label{e1}
([[a, b]\vartheta, [a, b]])^{m} = [[a, b]\vartheta, [a, b]],  ~\text{for all}  ~a, b \in \mathscr{Q}.
\end{equation}
By Martindale's theorem \cite{24}, it observes that $Soc(\mathscr{Q}) = F  \neq (0)$ and $eFe$ is finite  dimensional simple central algebra over $\mathscr{C}$, for any minimal idempotent element $e \in F$. We can also suppose that $F$ is non-commutative, because else $\mathscr{Q}$ must be commutative. Clearly, $F$ satisfies $([[a, b]\vartheta, [a, b]])^{m} = [[a, b]\vartheta, [a, b]]$ (see, for example, the proof of \cite[Theorem~1]{23}). As $F$ is a simple ring, either $F$ does not contain any non-trivial idempotent element or $F$ is generated by its idempotents. In this last case,
assume that $F$ contains two minimal orthogonal idempotent elements $e$ and $f$. Using the assumption, one can see that, for $[a,b] = [ea, f] = eaf$, we have
\begin{equation}\label{e2}
eaf\vartheta eaf =0,
\end{equation}
in this case we get $f\vartheta eaf \vartheta eaf\vartheta e = 0$, and  primeness of $\mathscr{R}$, we get $f\vartheta e = 0$ for any rank $1$ orthogonal idempotent element $e$ and $f$. Notably, for any rank $1$ idempotent element $e$, we have $e\vartheta (1 - e)=0$ and $(1 - e)\vartheta e=0$, that is, $e\vartheta=e\vartheta e=\vartheta e$. Hence, $[\vartheta,e]=0$
gives that $F$  is commutative or  $\vartheta \in \mathscr{C}$. We get a contradiction, in this case.

Now, we consider  when $F$ cannot contain two minimal orthogonal idempotent elements and so, $F = D$ for suitable finite dimensional division ring $D$ over its center which implies that $\mathscr{Q} = F$ and $\vartheta \in F$. By \cite[Theorem 2.3.29]{20} (see also \cite[Lemma 2]{23}), there exists a field  $\mathbb{K}$  such that $F \subseteq  M_{n}(\mathbb{K})$ and $M_{n}(\mathbb{K})$ satisfies
$([[a, b]\vartheta, [a, b]])^{m} = [[a, b]\vartheta, [a, b]]$. If $n = 1$ then $F\subseteq  \mathbb{K}$ and we have also a contradiction. Moreover, as we have just seen, if $n \geq  2$, then $\vartheta \in  Z(M_{n}(\mathbb{K}))$.

Finally, if $F$ does not contain any non-trivial idempotent element, then $F$ is finite dimensional division algebra over $\mathscr{C}$ and $\vartheta \in F  = \mathscr{R}\mathscr{C} = \mathscr{Q}$. If $\mathscr{C}$ is finite, then $F$ is finite division ring, that is, $F$ is a commutative field and so $\mathscr{R}$ is commutative too. If $\mathscr{C}$ is infinite, then $F \bigotimes_{\mathscr{C}} \mathbb{K} \cong M_{n}(\mathbb{K})$, where $\mathbb{K}$ is a splitting field of $F$. We  get the conclusion.

Henceforward, $\mathscr{T}$ is non-identity automorphism of $\mathscr{R}$.
Now, we have two cases:

Case I: If $\mathscr{T}$ is $\mathscr{Q}-$inner, then there exists an invertible element $\alpha$ of $\mathscr{Q}$ such that $\mathscr{T}(a) = \alpha a \alpha^{-1}$ for every $a \in \mathscr{R}$. Thus,  $([\alpha[a, b]\alpha^{-1}\vartheta, [a, b]])^{m} = [\alpha[a, b] \alpha^{-1}\vartheta, [a, b]]$  for every $a, b \in \mathscr{R}$.
If $\alpha^{-1}\vartheta \in \mathscr{C}$, then $\mathscr{R}$ satisfies $([\alpha[a, b], [a, b]])^{m} = [\vartheta[a, b], [a, b]]$ and we get the conclusion as above. Now we assume that $\alpha^{-1}\vartheta \not \in \mathscr{C}$, therefore
\[
([\alpha[a,b]\alpha^{-1}\vartheta, [a, b]])^{m} = [\alpha[a, b], [a, b]]
\]
 is a non-trivial GPI  for $\mathscr{R}$ and hence for $\mathscr{Q}$ by Fact \ref{f2}. In light of  ``Martindale's theorem \cite{24}, $\mathscr{Q}$ is isomorphic  to a dense subring of linear transformations of a vector space $V$ over $D$, where $D$ is a finite dimensional division ring over $\mathscr{C}$''. By Lemma \ref{l1},  the result follows.

Case II: If $\mathscr{T}$ is $\mathscr{Q}$-outer, and $\mathscr{Q}$ satisfies $([\mathscr{T}([a, b])\vartheta, [a, b]])^{m} = [\mathscr{T}([a, b])\vartheta, [a, b]]$, then by Lemma \ref{l1} we get $dim_{D}V = 1$, that is $\mathscr{Q}$ is a domain. By Fact \ref{f3}, $\mathscr{Q}$ satisfies $[[r, s]\vartheta, [a, b]]^{m} = [[r, s], [a, b]]$ and in particular, for $r=a$ and $s = b$, we have $[[a, b]\vartheta, [a, b]]^{m} = [[a, b]\vartheta, [a, b]]$  for every $a, b \in \mathscr{Q}$. Hence, using the same technique as above we get the  required  conclusion.
\end{proof}

\begin{theorem}\label{t1}
Let  $\mathscr{R}$ be a prime ring of Char$(\mathscr{R}) \neq 2$ and  $m\neq 1$ be a positive integer.  If $S$ is a nonzero skew derivation with an associated automorphism $\mathscr{T}$ of $\mathscr{R}$ such that $([S([a, b]), [a, b]])^{m} = [S([a, b]), [a, b]]$ for all $a, b \in \mathscr{R}$, then $\mathscr{R}$ is commutative.
\end{theorem}

\begin{proof}
We have
$$([S([a, b]), [a, b]])^{m} = [S([a, b]), [a, b]] ~\text{for every} a, b \in \mathscr{R}.$$
Firstly, we assume that $S$ is $\mathscr{Q}$-inner, that is, $S(a) = \vartheta a - \mathscr{T}(a)\vartheta$ with $0 \neq \vartheta \in \mathscr{Q}$. Thus, $\forall a, b \in \mathscr{R}$, we have
\[
[\vartheta[a, b] - \mathscr{T}([a, b])\vartheta, [a, b]])^{m} = [\vartheta[a, b] - \mathscr{T}([a, b])\vartheta, [a, b]].
\]
 If $\vartheta \in \mathscr{C}$, then $\mathscr{R}$ satisfies the GPI $([\mathscr{T}([a, b])\vartheta, [a, b]])^{m} = [\mathscr{T}([a, b])\vartheta, [a, b]]$.  We get the desired conclusion,  by Proposition \ref{p1}. Therefore $\vartheta \not \in \mathscr{C}$, and so \[
 [\vartheta[a, b] - \mathscr{T}([a, b])\vartheta, [a, b]])^{m} = [\vartheta[a, b] - \mathscr{T}([a, b])\vartheta, [a, b]]
 \]
 is nontrivial GPI for $\mathscr{R}$. Thus, Lemma \ref{l1} yields the required result.

Finally,  when $S$ is $\mathscr{Q}$-outer, then the above identity can be rewritten as
\[
[S(a)b + \mathscr{T}(a)S(b) - S(b)a \mathscr{T}(b)S(a), [a, b]]^{m} = [S(a)b + \mathscr{T}(a)S(b) - S(b)a -\mathscr{T}(b)S(a), [a, b]],
\]
 and hence $\mathscr{R}$ satisfies
\[
([\vartheta b+\mathscr{T}(a)s-sa-\mathscr{T}(b)r, [a, b]])^{m} = [rb+\mathscr{T}(a)s -sa -\mathscr{T}(b)r, [a, b]].
\]
 In particular $\mathscr{R}$ satisfies $([\mathscr{T}(a)s-sa, [a, b]])^{m} = [\mathscr{T}(a)s -sa, [a, b]]$.  We divide it into  two cases.
First, $\mathscr{T}$ be an identity map of $\mathscr{R}$.  Then $([[r,s],[a,b]])^{m} = [[r,s], [a,b]]$ for every $a, b, r, s \in \mathscr{R}$, that is, $\mathscr{R}$ is a polynomial identity ring. Thus, $\mathscr{R}$ and $M_{n}(\mathbb{K})$ satisfy the same polynomial identities \cite[Lemma 1]{23}, i.e.,
\[
([[r,s],[a,b]])^{m} = [[r,s], [a, b]]
\quad \textup{ for each }
a,b,r,s \in M_{n}(\mathbb{K}),
\]
 Let $n\geq 2$ and $e_{ij}$ be the usual unit matrix. Then $r=b=e_{12}$, $s=e_{21}$ and $a=e_{11}$, we get a contradiction $2e_{12} =0$. Thus, $n =1$ and we are done.

Now consider $\mathscr{T}$  is not the identity map. Therefore,
\[
([\mathscr{T}(a)s - sa,[a, b]])^{m} = [\mathscr{T}(a)s - sa, [a, b]]
\]
 is a non-trivial GPI for $\mathscr{R}$, by Main Theorem in \cite{5}. Moreover, by Fact \ref{f2},
$\mathscr{R}$ and $\mathscr{Q}$ satisfy the same GPIs with automorphisms and
hence $([\mathscr{T} (a)s  - sa,[a, b]])^{m} = [\mathscr{T}(a)s -sa, [a,b]]$ is also an identity for $\mathscr{Q}$. Since $\mathscr{R}$ is a GPI-ring, by \cite{24}
``$\mathscr{Q}$ is a primitive ring, which is isomorphic to a dense subring of the ring of linear
transformations of a vector space $V$ over a division ring $D$''.
If  $\mathscr{Q}$ is a domain, then by Fact \ref{f3}, we have that $\mathscr{Q}$  satisfies the equation $([ts - sa, [a, b]])^{m} =[ts - sa, [a, b]]$.  In
particular, $([[a,z],[a, b]])^{m}  = [[a, z], [a, b]]$ for all $a, b,z \in  \mathscr{Q}$, which yields that
$\mathscr{Q}$ is commutative (by using the same above argument) and hence $\mathscr{R}$.  Henceforth, $\mathscr{Q}$ is not a domain. We have $\mathscr{T}(a) = hah^{-1}$   $\forall a \in \mathscr{Q}$,  as mentioned above. Thus,  $([hah^{-1}z - za, [a, b]])^{m} = [hah^{-1}z - za, [a, b]]$
Hence, we may consider  that
$dim~D_V \geq  2$.
By  \cite[p. 79]{20}, there exists a semi-linear automorphism $h \in End(V)$ such that
$\mathscr{T}(a) = hah^{-1}$  $\forall a \in \mathscr{Q}$. Hence, $\mathscr{Q}$ satisfies  $([hah^{-1}z - za, [a, b]])^{m} = [ hah^{-1}z -za, [a, b]]$.

If for any $v \in V$ $\exists ~\Theta_v \in D$ such that $h^{-1}v = v \Theta_{v}$, then, it follows that there exists a unique $\Theta  \in D$ such that $h^{-1}v = v\Theta$,
$\forall v \in V$ (see for example Lemma 1 in \cite{ccl}). In this case
$\mathscr{T}(a)v = (hah^{-1})v = hav\Theta$
and
$$(\mathscr{T}(a) - a)v = h(av\Theta) - av = h(h^{-1}av) - av = 0,$$
since $V$ is faithful, which  is a  contradiction that $\mathscr{T}$  is the identity map.
Thus,  $\exists$ $v \in V$ such that $\{v, h^{-1}v\}$ is linearly $D$-independent.
In this case, first we assume that  $dim ~V_{D} \geq 3$. Thus, $\exists~ u \in V$ such that $\{u, v,  h^{-1}v \}$ is linearly $D$-independent. Hence, the density theorem for $\mathscr{Q}$,  $\exists ~a, b, z \in \mathscr{Q}$ such that
\begin{align*}
zv &= 0 & zh^{-1}v &= h^{-1}v\\
bv &=0 & bh^{-1}v &= 0\\
av &= h^{-1}v & bu &= -2v\\
ah^{-1}v &= u .&&
\end{align*}
The above relation gives that
$$0 = (([hah^{-1}z -za, [a, b]])^{m} - [hah^{-1}z - za, [a, b]])v = (2^{m} - 2) v \neq 0$$
again a contradiction.

Now, the case when $dim~V_{D} = 2$ that is, $\mathscr{Q} = M_{2}(\mathbb{K})$. Thus
\[
([\mathscr{T}(a)z -za, [a, b]])^{2} =[\mathscr{T}(a)z -za, [a, b]]
\quad \textup{ for all } a, b, z \in \mathscr{Q}.
\]
 Since $\mathscr{T}(a)$-word of degree  $2$  and Char$(\mathscr{R}) >3$ by \cite[Theorem 3]{6}, 
 \[
 ([tz-za, [a,b]])^{2} - [tz -za, [a, b]] = 0
 \quad \textup{ for every } t, z, a, b \in \mathscr{Q}.
 \]
  Using the same technique as above  its shows that $\mathscr{Q}$ is commutative and hence $\mathscr{R}$ is commutative.
\end{proof}

The following corollary  is an immediate consequence of our result.
\begin{corollary}
	\cite[Theorem 2.3]{10}
	Let $\mathscr{R}$ be a prime ring of characteristic not two and $d$ be a nonzero derivation of $\mathscr{R}$  such that
	$ ([d([a, b]), [a, b]])^{m} = [d([a, b]), [a, b]]$ for all $a, b \in  \mathscr{R}$. Then $\mathscr{R}$ is commutative.
\end{corollary}

\begin{theorem}\label{t2}
	Let $\mathscr{R}$ be a prime ring of Char$(\mathscr{R}) \neq 2$, $m\neq 1$ be a positive integer   and $\mathscr{L}$ a Lie ideal of $\mathscr{R}$. If $S$ is a nonzero skew derivation with an associated automorphism $\mathscr{T}$ of $\mathscr{R}$ such that $([S(v), v])^{m} = [S(v), v]$ for all $v \in \mathscr{L}$, then $L$ contained in the center of $\mathscr{R}$.
\end{theorem}

\begin{proof}
 Suppose that $\mathscr{L} \not \subseteq  Z(\mathscr{R})$ is a Lie ideal of $\mathscr{R}$. Then by \cite{15}, there exists an ideal $I$ of $\mathscr{R}$ such that $0 \neq [I, \mathscr{R}] \subseteq \mathscr{L}$ and $[\mathscr{L}, \mathscr{L}] \neq (0)$.
 Also, $\mathscr{R} \not \subseteq Z(\mathscr{R})$ as $\mathscr{L}$ is a noncentral Lie ideal of $\mathscr{R}$. Therefore  by the given hypothesis, $I$ as well as $\mathscr{R}$ (Fact \ref{f2}) satisfy $[S([a, b]), [a, b]])^{m} = [S([a, b]), [a, b]]$. By Theorem \ref{t1}, we get the required result.
\end{proof}

\subsection*{Acknowledgment:}
The authors are greatly indebted to the referee for his/her valuable suggestions, which have immensely improved the paper. For the first author, this research is supported by the Council of Scientific and Industrial Research (CSIR-HRDG), India, Grant No.  25(0306)/20/EMR-II.


\EditInfo{June 17, 2021}{August 12, 2021}{Ivan Kaygorodov}

\end{paper}